\documentclass[graybox]{svmult}

\usepackage{mathptmx}       
\usepackage{helvet}         
\usepackage{courier}        
\usepackage{subfigure}
\usepackage{bm}
\usepackage{makeidx}         
\usepackage{graphicx}        
\usepackage{multicol}        
\usepackage[bottom]{footmisc}
\usepackage{amsmath,amsfonts,amssymb}
\usepackage{subeqnarray}
\usepackage{color}
\usepackage[usenames,dvipsnames]{xcolor}

\usepackage[numbers]{natbib}

\makeindex             
                       
  \newcommand\solidrule[1][5mm]{\rule[0.5ex]{#1}{2pt}}

  \newcommand{\bcal}[1]{{\cal{#1}}}

  \newcommand{\be}{\begin{equation}}
    \newcommand{\ee}{\end{equation}}
  \newcommand{\bea}{\begin{eqnarray}}
    \newcommand{\eea}{\end{eqnarray}}
    \newcommand{\beas}{\begin{subeqnarray}}
    \newcommand{\eeas}{\end{subeqnarray}}
  \newcommand{\ube}{\begin{equation*}}
    \newcommand{\uee}{\end{equation*}}
  \newcommand{\ubea}{\begin{eqnarray*}}
    \newcommand{\ueea}{\end{eqnarray*}}
  \newcommand{\ben}{\begin{enumerate}}
    \newcommand{\een}{\end{enumerate}}
  \newcommand{\bit}{\begin{itemize}}
    \newcommand{\eit}{\end{itemize}}

  \newcommand{\ext}{\mbox{d}}
  \newcommand{\hodge}{\star}

  \newcommand{\coD}{d^\hodge}
  
  \newcommand{\bounda}{\partial}

  \newcommand{\partiald}[2]{\dfrac{\partial #1}{\partial #2}}
  
  \newcommand{\intP}[1]{{\mathfrak{i}}_{\vec{#1}}}

  \newcommand{\innerP}[2]{\left ( #1 , #2 \right)_{\phyD}}

  \newcommand{\phyD}{\Omega}

  \newcommand{\vB}[1]{{\bm{\partial_{#1}}}}

  \newcommand{\figref}[1]{{Fig.~\ref{#1}}}

\begin{document}

\title*{A Geometric Approach Towards Momentum Conservation}
\author{D. Toshniwal, R.H.M. Huijsmans and M.I. Gerritsma}
\authorrunning{D. Toshniwal, R. Huijsmans and M. Gerritsma}
\institute{Deepesh Toshniwal \at \;
\email{toshniwald.iitkgp@gmail.com}
\and R.H.M. Huijsmans \at Maritime Engineering, TU Delft, The Netherlands.
\email{R.H.M.Huijsmans@tudelft.nl}
\and M.I. Gerritsma \at Aerospaee Engineering, TU Delft, The Netherlands.
\email{M.I.Gerritsma@tudelft.nl}
}
%
%
\maketitle

\abstract*{In this work, a geometric discretization of the
Navier-Stokes equations is sought by treating momentum as a
covector-valued volume-form. The novelty of this approach is that we treat conservation of momentum as a tensor equation and describe a higher order approximation to this tensor equation. The resulting scheme satisfies mass and momentum
conservation laws exactly, and resembles a staggered-mesh finite-volume method.
Numerical test-cases to which the discretization scheme is applied are the
Kovasznay flow, and lid-driven cavity flow.}

\abstract{In this work, a geometric discretization of the
Navier-Stokes equations is sought by treating momentum as a
covector-valued volume-form. The novelty of this approach is that we treat conservation of momentum as a tensor equation and describe a higher order approximation to this tensor equation. The resulting scheme satisfies mass and momentum
conservation laws exactly, and resembles a staggered-mesh finite-volume method.
Numerical test-cases to which the discretization scheme is applied are the
Kovasznay flow, and lid-driven cavity flow.}

\section{Navier-Stokes equations} \label{sec:navierStokes}
Mimetic discretizations aim to represent physics in a discrete sense, in contrast to differential formulations, which are concerned with the limit $h \rightarrow 0$. For the case in which $h \neq 0$ geometrical considerations play an important role in the correct discrete formulation, \citep{Bossavit1,Bochev1,Desbrun1,marcI,KreeftF}. Application of these ideas to continuum models are described in \citep[Appendix A]{Frankel} and \citep{Kanso, Yavari}. The novel aspect in this paper is that continuum ideas are applied to incompressible, viscous flows using spectral basis functions.
%

We start with the incompressible Navier-Stokes equations ($\rho=1$), written in the integral formulation, as given in many textbooks and we try to make precise what these statements mean. It is important to give an accurate meaning to all variables, because when we want to represent these physical quantities on finite grids, we want to preserve the main structure of the equations.
Conservation of mass ($\rho=1$) is usually given by
\be
  \int_{\bounda \phyD} \vec{v} \cdot \vec{n}\,\mathrm{d}S = 0\;,
\label{eq:NavierStokes-Mass}
\ee
and conservation of momentum,
\be
   \int_{\phyD}\partiald{ \vec{v}}{t}\;dV + \int_{\bounda \phyD}\vec{v} \otimes \vec{v} \cdot \vec{n}\;dS =
\int_{\bounda \phyD} \mathbf{\sigma} \cdot \vec{n}\;dS\;
\label{eq:NavierStokes-Momentum}
\ee
and Newtonian stress relation
\be
\mathbf{\sigma} = -p \mathbb{I} + \mu \left ( \nabla \vec{v} + \left ( \nabla \vec{v} \right )^T \right ) \;.
\ee
Here $\vec{v}, p, \sigma \mbox{ and } \mu $ denote velocity, pressure, total stress tensor and dynamic viscosity,
respectively; $\mathbb{I}$ and $\vec{n}$ are the identity matrix and the outward unit normal to the boundary, respectively.
The above are balance equations for volumetric quantities that depend on their fluxes through surfaces and are more physical than their differential counterparts. 
\subsection{Momentum and velocity}
The first term in (\ref{eq:NavierStokes-Momentum}) indicates that velocity (and its time derivative) can be integrated over a volume. But velocity is generally {\em not} associated to volumes, but is defined as the tangent vector at a given point along the trajectory of a particle. Velocity is therefore a vector-valued $0$-form. This statement means that to every point in space-time (a zero-dimensional object) we associate a vector. Let $V$ be the linear vector-space of all possible vectors at a given point in space, then we can define the space $V^*$ of all linear functionals on $V$. Elements of $V^*$ are called {\em covectors}. The spaces $V$ and $V^*$ are isomorphic, but there is no canonical isomorphism which relates an element $v \in V$ to an element $\alpha \in V^*$. Once a metric is defined, one can associate with every vector at a point a corresponding covector. This map is called the flat operator: $\flat \,: V \rightarrow V^*$. The covector associated with a vector $v$ is then denoted by $v^\flat$. 

The linear vector space $V$ associated to a point $p$ is called the {\em tangent space} at $p$, denoted by $T_p\phyD$. The corresponding dual space is called the {\em cotangent space} at $p$ denoted by $T_p^*\phyD$. The collection of all tangent spaces in the domain $\phyD$ is called the {\em tangent bundle}, $T\phyD$ and the collection of cotangent spaces is called the {\em cotangent bundle}, $T^*\phyD$. Let $\alpha \in T^*\phyD$ and $\vec{v} \in T\phyD$, then $\langle \alpha, \vec{v} \rangle$ associates to each point $p$ in $\phyD$ the value $\left. \alpha \right |_p(\left. \vec{v} \right |_p )$.

With every $k$-form we can associate a $(n-k)$-form with a different type of orientation, see \citep{Bossavit1}. The collection of all $k$-forms on $\phyD$ is denoted by $\Lambda^k(\phyD)$. The metric dependent operator which establishes this connection is the Hodge-$\star$ operator. For continuum models we need to combine the $\flat$ and Hodge-$\star$ into the operator $\star^\flat$, (see also \citep{Yavari} for such operations)
\[ \star^\flat \, : \, T\phyD \otimes \Lambda^k(\phyD) \rightarrow T^*\phyD \otimes \Lambda^{n-k}(\phyD) \;.\]
If we apply this operator to velocity $\vec{v} \in T\phyD \otimes \Lambda^0(\phyD)$ we obtain
\[ m:=\star^\flat(\vec{v}) \in T^*\phyD \otimes \Lambda^n{\phyD} \;.\]
Similarly, we can define $\star^\sharp \, : \, T^*\phyD \otimes \Lambda^k(\phyD) \rightarrow T\phyD \otimes \Lambda^{n-k}(\phyD)$.
The physical quantity $m$ is called {\em momentum density} or the {\em momentum per unit volume}. This is a {\em covector-valued volume form}. So instead of integrating `velocity' over the domain
we are tempted to write
\[ \int_{\phyD} m = \int_\phyD \star^\flat(\vec{v}) \;.\]
This integral is not defined, because it assumes that we can integrate over the tangent spaces in $\phyD$. The basis in each tangent space, however, may differ from point to point. In order to define the momentum integral we introduce the operator $\stackrel{\cdot}{\wedge}$
\[ \stackrel{\cdot}{\wedge} \, :\, \left ( T^*\phyD \otimes \Lambda^k(\phyD) \right ) \otimes \left ( T\phyD \otimes \Lambda^l(\phyD) \right ) \rightarrow \Lambda^{k+l}(\phyD)\;,\]
given by $\alpha \in T^*\phyD \otimes \Lambda^k(\phyD)$ and $\vec{w} \in T\phyD \otimes \Lambda^l(\phyD)$
\[ \alpha \stackrel{\cdot}{\wedge} \vec{w} = \langle \alpha,\vec{w} \rangle \mathrm{d}x^{(k)}\wedge \mathrm{d}x^{(l)} \;.\]
This operation yields a $(k+l)$-form which can be integrated over $(k+l)$-dimensional submanifolds.

If we apply momentum density $m$ to any vector field $\vec{w}$ (not necessarily a velocity field) using this operator we get $m \stackrel{\cdot}{\wedge} \vec{w} \in \Lambda^n(\phyD)$ and this can be integrated over a volume. So the proper way to interpret the time rate of change of momentum should be
\begin{equation}
 \int_\phyD \frac{\partial}{\partial t} \star^\flat(\vec{v}) \stackrel{\cdot}{\wedge} \vec{w} \;,\;\;\; \forall \vec{w} \in T\phyD \otimes \Lambda^0(\phyD)\;.
 \label{eq:correct_momentum_change}
 \end{equation}
In many textbooks on fluid the distinction between dynamics momentum density (usually called `momentum') and velocity is ignored; one is just a scalar multiple of the other, $m = \rho \vec{v}$, but the use of the vector $\vec{w}$ in (\ref{eq:correct_momentum_change}) is generally incorporated. The textbooks then say: `We consider this equation for each component separately ..'. This is a strange sentence, because components have no physical relevance, only vectors, i.e. components plus associated basis vectors are physically relevant. But what is meant by this statement is that for the vector field $\vec{w}$ in (\ref{eq:correct_momentum_change}) a uniform vector field in the $x^i$-direction is taken. The generality `all vector fields' is in these textbooks compensated by the fact that momentum conservation should hold for `all volumes'.

\subsection{Convection}
Now that we understand how momentum density should be integrated over a volume, we can also define convection of momentum density. After pairing with an arbitrary vector field, $\vec{w}$, we obtain a volume form and we apply the Lie derivative to this volume form, see \citep{palhaI}. The Lie derivative for a volume form, $\beta^{(n)}$, is given by
\[ {\mathcal L}_{\vec{v}} \beta^{(n)} = \mathrm{d} \intP{v} \beta^{(n)} \;,\]
and then the generalized Stokes theorem converts this exact form to a boundary integral
\begin{equation} 
 \int_{\phyD} {\mathcal L}_{\vec{v}} m \stackrel{\cdot}{\wedge} \vec{w} = \int_{\partial \phyD} \intP{v} ( m \stackrel{\cdot}{\wedge} \vec{w}) \;.\label{eq::convMomRelation}
\end{equation}
Compare this expression with the the convective term in (\ref{eq:NavierStokes-Momentum}) and note that it does not require an inner product nor the definition of an outward unit normal. The inner product is avoided since we work with differential forms and duality pairing is metric-free and the orientation of the elements in the mesh, \citep{marcI}, avoids the use of explicitly defined normals.

\subsection{Stress tensor and surface force density}
The last term in (\ref{eq:NavierStokes-Momentum}) denotes the action of the viscous forces on the flow represented by the stress tensor $\mathbf{\sigma}$. The stress tensor is an infinitesimal quantity in the limit for $h \rightarrow 0$. On a finite mesh we can identify volumes over which we integrate the momentum density and the boundary of these volumes where surface forces act. In continuum mechanics forces are `smeared out', so we introduce the {\em surface force density} given by $\mathfrak{t} \in T^*\phyD \otimes \Lambda^{n-1}(\phyD)$. This is a {\em covector-valued $(n-1)$-form}. Forces are generally associated with covectors, \citep{Bossavit1,Tonti1}, and in the current setting need to be covectors in order to equate them to the time rate of change of momentum which was also covector-valued. It is furthermore a $(n-1)$-form since it acts on the the boundary of $n$-dimensional volumes, see also \citep[Appendix A]{Frankel} and \citep{Kanso,Yavari}. Again, covector-valued forms cannot be integrated, so the 
proper way is pair this to pair it with an arbitrary vector field $\vec{w}$ before integration over surfaces is possible. The momentum equation then becomes
\be
\frac{d}{dt} \int_\phyD \star^\flat(v) \stackrel{\cdot}{\wedge} \vec{w} + \int_{\partial \phyD} \intP{v} ( \star^\flat(v) \stackrel{\cdot}{\wedge} \vec{w}) = \int_{\partial \phyD} \mathfrak{t} \stackrel{\cdot}{\wedge} \vec{w} \;,\;\;\; \forall \vec{w} \in T\phyD \otimes \Lambda^0(\phyD)\;.
\ee

\subsection{Newtonian stress relation}
The pressure scalar is an outer-oriented volume form, $p^{(n)}$.
Pressure force density is represented as a covector-valued $(n-1)$-form
\[ \vec{p} = \left ( \star p \right ) \, \mathrm{d}x^i \otimes \mathrm{d}x^1 \wedge \ldots \widehat{\mathrm{d}x^i} \ldots \wedge \mathrm{d}x^n\;,\]
where the notation $\widehat{\cdot}$ indicates that this term is omitted and $\mathrm{d}x^i \otimes \mathrm{d}x^1 \wedge \ldots \widehat{\mathrm{d}x^i} \ldots \wedge \mathrm{d}x^n$ is the identity tensor, see also example \citep[\S 9.3a]{Frankel}. This description agrees with \citep{KreeftJCP} for Stokes flow. Note that $\vec{p} \stackrel{\cdot}{\wedge} \vec{w} = \intP{w} p^{(n)}$.

The velocity gradient is represented as the covariant differential of the velocity vector field, $\nabla \vec{v}$ which is a vector-valued $1$-form, see \citep[\S 9.3b]{Frankel}. In this paper we restrict ourselves to Euclidean space for which the connection $1$-forms vanish. Applying $\star^\flat_\mu (\nabla \vec{v})$ transforms the vector-valued $1$-form into a covector-valued $(n-1)$-form, where the diffusion coefficient is contained in the Hodge-$\star$ operator. In this paper we assume $\mu$ to be constant.

\subsection{Conservation of mass}
Let $\omega^{(n)}$ be the standard volume form, then the divergence of a vector field is defined as $(\mbox{div}\,\vec{v}) \omega^{(n)}= {\mathcal L}_{\vec{v}} \omega^{(n)} = \mathrm{d} \intP{v} \omega^{(n)}$. Integration over a volume and applying Stokes theorem gives
\[ \int_\phyD \mathrm{d} \intP{v} \omega^{(n)} = \int_{\partial \phyD} \intP{v} \omega^{(n)} \;.\]
This is the proper translation of (\ref{eq:NavierStokes-Mass}) as found in textbooks on incompressible flow. The velocity flux field, $\intP{v} \omega^{(n)}$, is isomorphic to the velocity vector field. The velocity flux field will be used in the discrete representation of the Navier-Stokes equations. The relation between the velocity fluxes $\intP{v} \omega^{(n)}$ and $\star^\flat(\vec{v})$ is given by
\be
\star^\flat(\vec{v}) \stackrel{\cdot}{\wedge} \vec{w} = \vec{w}^\flat \wedge \intP{v} \omega^{(n)}\;,\;\; \forall \vec{w} \in T\phyD \otimes \Lambda^0(\phyD)\;.
\label{eq:momentum_velocity_flux}
\ee

%
%
%
The volume forms and $(n-1)$-forms appearing in all integrals are all {\em outer-oriented}.

\section{Discrete representation}
In the full differential geometric setting as described above, the integration only makes sense when paired with all vector fields $\vec{w}$. Here we choose the uniform vector field in the $x$- and $y$-direction only and impose that conservation should hold for all volumes in our spectral elements. These volumes are generated by the Gauss-Lobatto grid in the spectral element and will be denoted by $\Omega_{ij}$. So in this section $\vec{w}$ is either $\vB{x}$ or $\vB{y}$.
 


\begin{figure}[htb!]
\centering
  \includegraphics[trim=11mm 22mm 23mm 10mm,clip,width=0.85\linewidth]{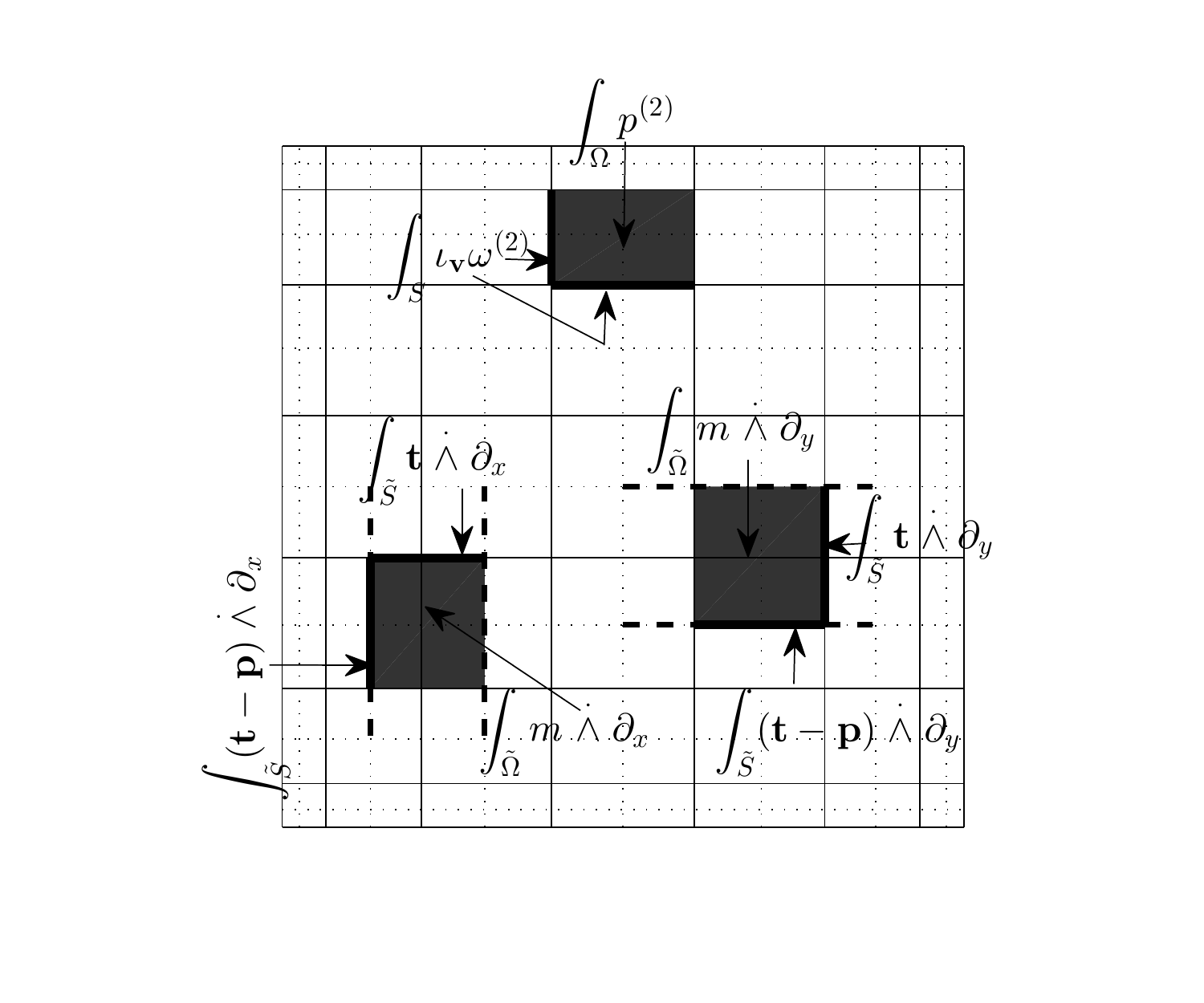}
\caption[Velocity and pressure placement on a lowest-order mesh]{The velocities are discretized as outer-oriented mass-fluxes, and live on surfaces ($S$) of the Gauss-Lobatto grid shown above, while pressure is discretized on volumes ($\Omega$). Momenta are discretized on staggered volumes ($\tilde{\Omega}$) and their fluxes on the surfaces  ($\tilde{S}$) surrounding these staggered volumes.}
\label{fig::discreteRep}
\end{figure}

Figure~\ref{fig::discreteRep} displays one spectral element and its Gauss-Lobatto grid (solid lines) in 2D. The dotted gray lines represent the dual grid, see \citep{marcI,KreeftF}.

Momentum is reduced onto a volume consisting of a primal $(n-1)$-chain and a dual $1$-chain. In 2D these volumes consist of tensor products of primal and dual edge as shown in Figure~\ref{fig::discreteRep} by volumes enclosed by solid (primal) and dashed (dual) lines. The location of the unknowns coincides with those in staggered finite volume methods. The difference is that in this formulation the unknowns represent integral values, whereas in finite volume methods the unknowns either represent average or nodal values. Let us denote the primal surfaces by $S_i$, then discrete velocity is given by
\[ \bar{v}_i = \int_{S_i} \intP{\vec{v}} \omega^{(n)} \;.\]
This yields a metric-free description of conservation as mass as shown in \citep{KreeftF,GerritsmaEd}. The reduction of the pressure field is on outer-oriented volumes, see also \citep{KreeftF}. 

Integrals of momentum flux, $\bcal{F}^{(n-1)}$, pressure force, $\intP{\vec{w}}p^{(n)}$, and velocity gradients are represented on the boundary of the momentum volumes indicated in Figure~\ref{fig::discreteRep}.

Once we have the discrete variables for mass flux, momentum and pressure, we use the spectral element functions described in \citep{KreeftF,GerritsmaEd}, to interpolate these values in such a way that the integral values are preserved.

Using (\ref{eq:momentum_velocity_flux}) we can write the relation between momentum and velocity flux as
\be
\begin{split}
  \int_{\tilde{\Omega}_{ij}} m \stackrel{\cdot}{\wedge} \vec{w} - \int_{\tilde{\Omega}_{ij}} \vec{w}^\flat \wedge \intP{\vec{v}} \omega^{(n)}  = 0\;\longrightarrow
 \bar{m}_{\vec{w}} - P^m_{\vec{w}} \bar{u} = 0\;,
\end{split}
\ee
where $\tilde{\Omega}_{ij}$ are the volume where momentum is reduced, see \figref{fig::discreteRep}, and $P^m_{\vec{w}}$ is the matrix which maps discrete velocity (which is discretized as mass-fluxes) to discrete momentum (on the staggered-grid).  The discrete representation of momentum-flux, pressure force $\intP{\vec{w}}p^{(n)}$ and traction forces, $\star^\flat_\mu (\nabla_{\vec{w}}\vec{v})$ (which can be equivalently written as $\star^\flat_\mu (\nabla\star^\sharp\star^\flat\vec{v})\stackrel{\cdot}{\wedge} \vec{w} = \star^\flat_\mu (\nabla\star^\sharp m)\stackrel{\cdot}{\wedge} \vec{w}$, which, in Cartesian coordinates and with constant $\vec{w}$ becomes $\coD_\mu(m \stackrel{\cdot}{\wedge}\vec{w})$), are given by
\begin{itemize}
 \item Convective-flux, see \citep{palhaI}, $\bcal{F}^{(1)}_{\vec{w}} = \intP{v} ( m \stackrel{\cdot}{\wedge} \vec{w})$:
 \be
  \begin{split}
    \innerP{\bcal{F}^{(1)}_{\vec{w}}}{\beta^{(1)}}& - \innerP{m \stackrel{\cdot}{\wedge} \vec{w}}{\vec{v}^{\flat} \wedge \beta^{(1)}}\;,
    = 0\;,\\
    \longrightarrow &\tilde{M}_{11}\bar{\bcal{F}}_{\vec{w}} - \tilde{C}_{\vec{v}}\bar{m}_{\vec{w}} = 0\;.
  \end{split}
  \ee
  \item Pressure-force, $\bcal{H}^{(1)}_{\vec{w}} = \intP{w} p^{(2)}$:
  \be
    \begin{split}
      \bcal{H}^{(1)}_{\vec{w}} - p^{(2)}(\vec{w}) = 0
      \longrightarrow &\bar{\bcal{H}}_{\vec{w}} - P^p_{\vec{w}} \bar{p} = \tilde{B}_P\;.\\
    \end{split}
  \ee
 \item Diffusive-fluxes, $\bcal{T}^{(1)}_{\vec{w}} = \coD_\mu(m \stackrel{\cdot}{\wedge}\vec{w})$:
  \be
    \begin{split}
    \innerP{\bcal{T}^{(1)}_{\vec{w}}}{\beta^{(1)}} - & \innerP{m \stackrel{\cdot}{\wedge} \vec{w}}{\ext \beta^{(1)}} = -\int_{\partial \phyD} \beta^{(1)} \wedge \hodge (m \stackrel{\cdot}{\wedge} \vec{w}) \;, \\
    \longrightarrow &\tilde{M}_{11}\bar{\bcal{T}}_{\vec{w}} - \tilde{D}_{21}^T\tilde{M}_{22}\bar{m} = \tilde{B}_T\;.
  \end{split}
  \ee
\end{itemize}
The discrete continuity equation is given by
\be
  D_{21}\bar{u} = 0\;.
\ee
\begin{figure}[htb!]
\centering
  \subfigure[{Pressure, $h$-refinement}]{\includegraphics[trim = 40mm 90mm 40mm 100mm,clip,width=0.49\linewidth]{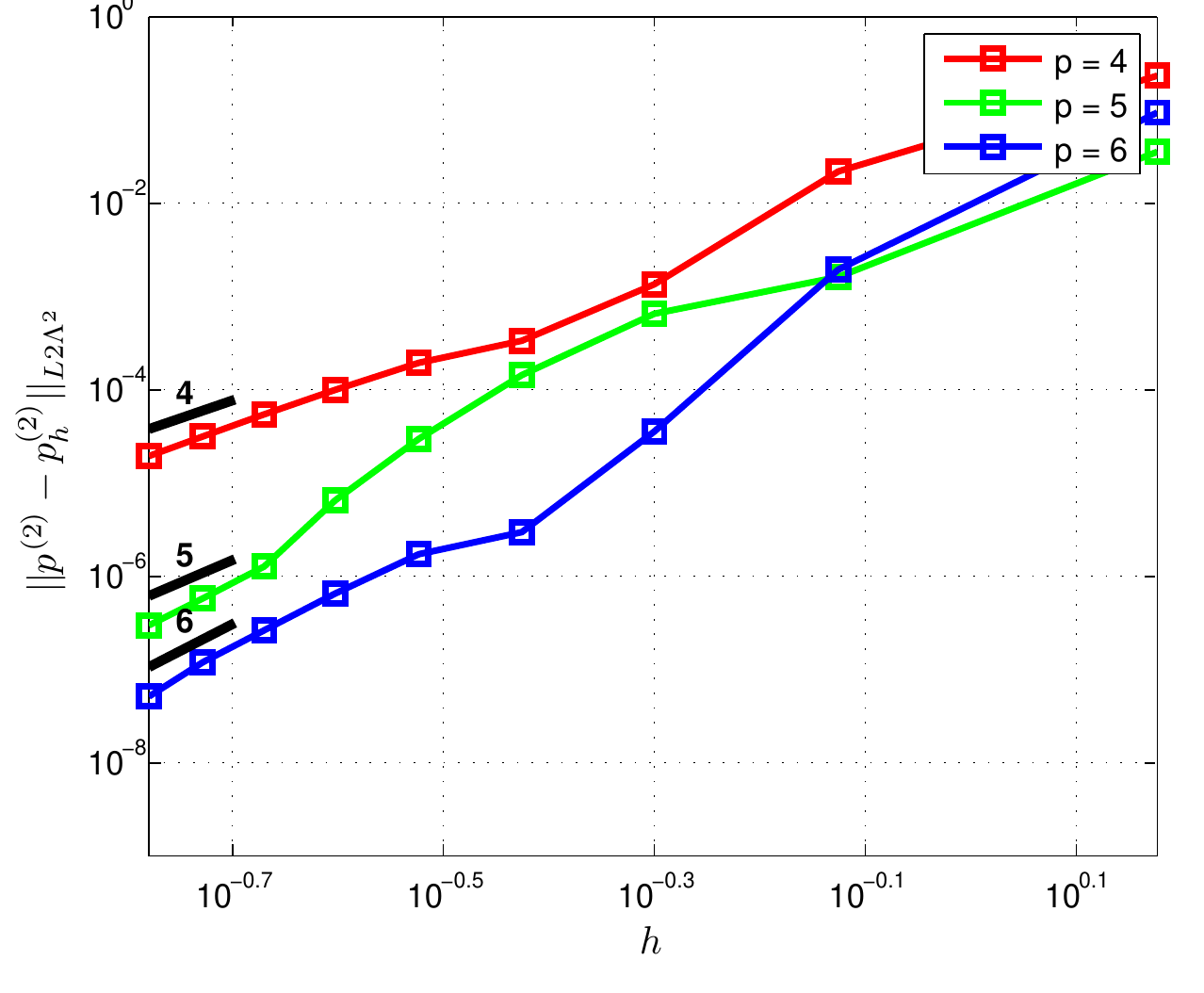}\label{fig::kovConva}} 
  \subfigure[{Pressure, $p$-refinement}]{\includegraphics[trim = 10mm 4mm 14mm 9mm,clip,width=0.49\linewidth]{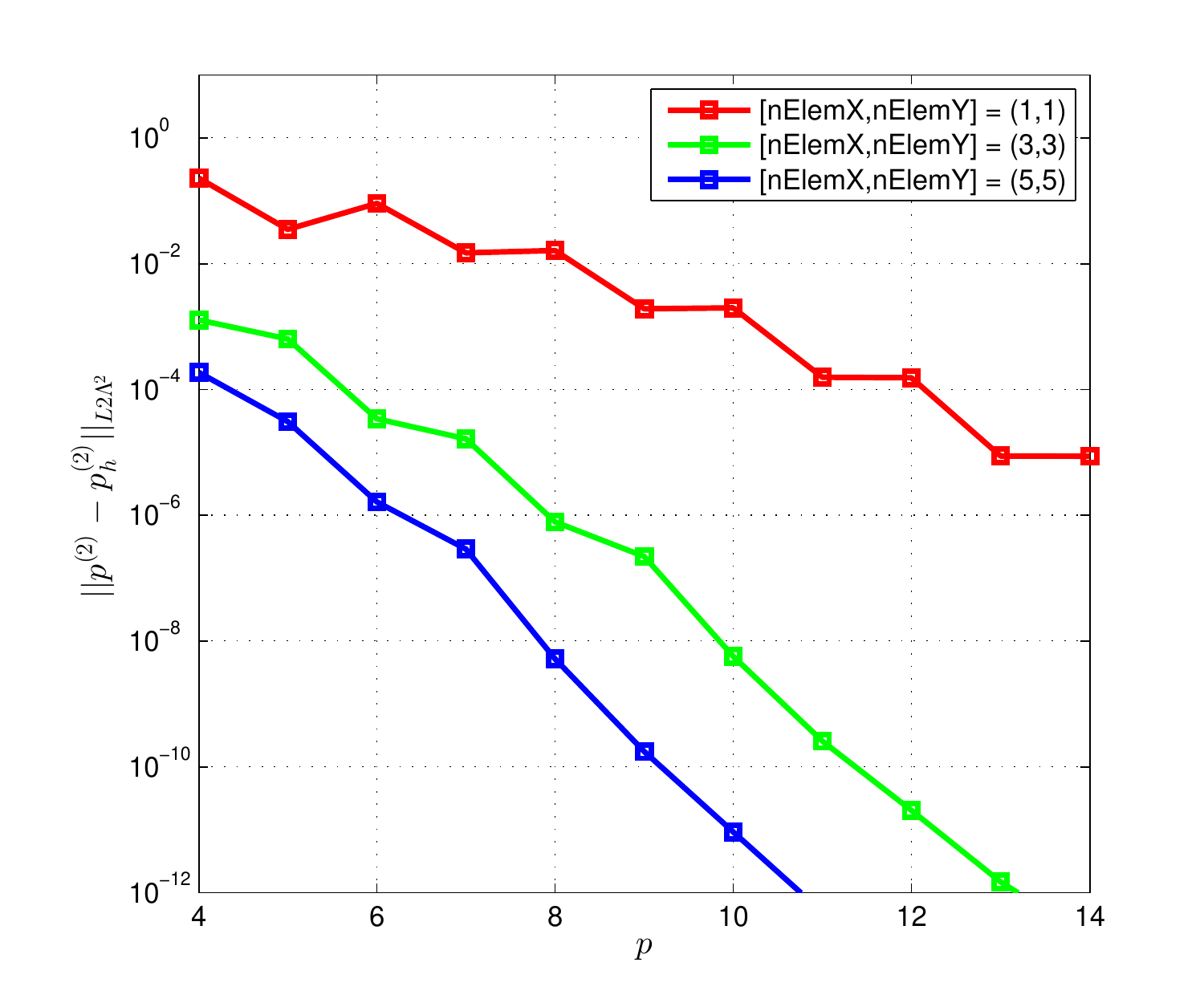}\label{fig::kovConvb}} \\
  \subfigure[{Velocity, $h$-refinement}]{\includegraphics[trim = 40mm 90mm 40mm 100mm,clip,width=0.49\linewidth]{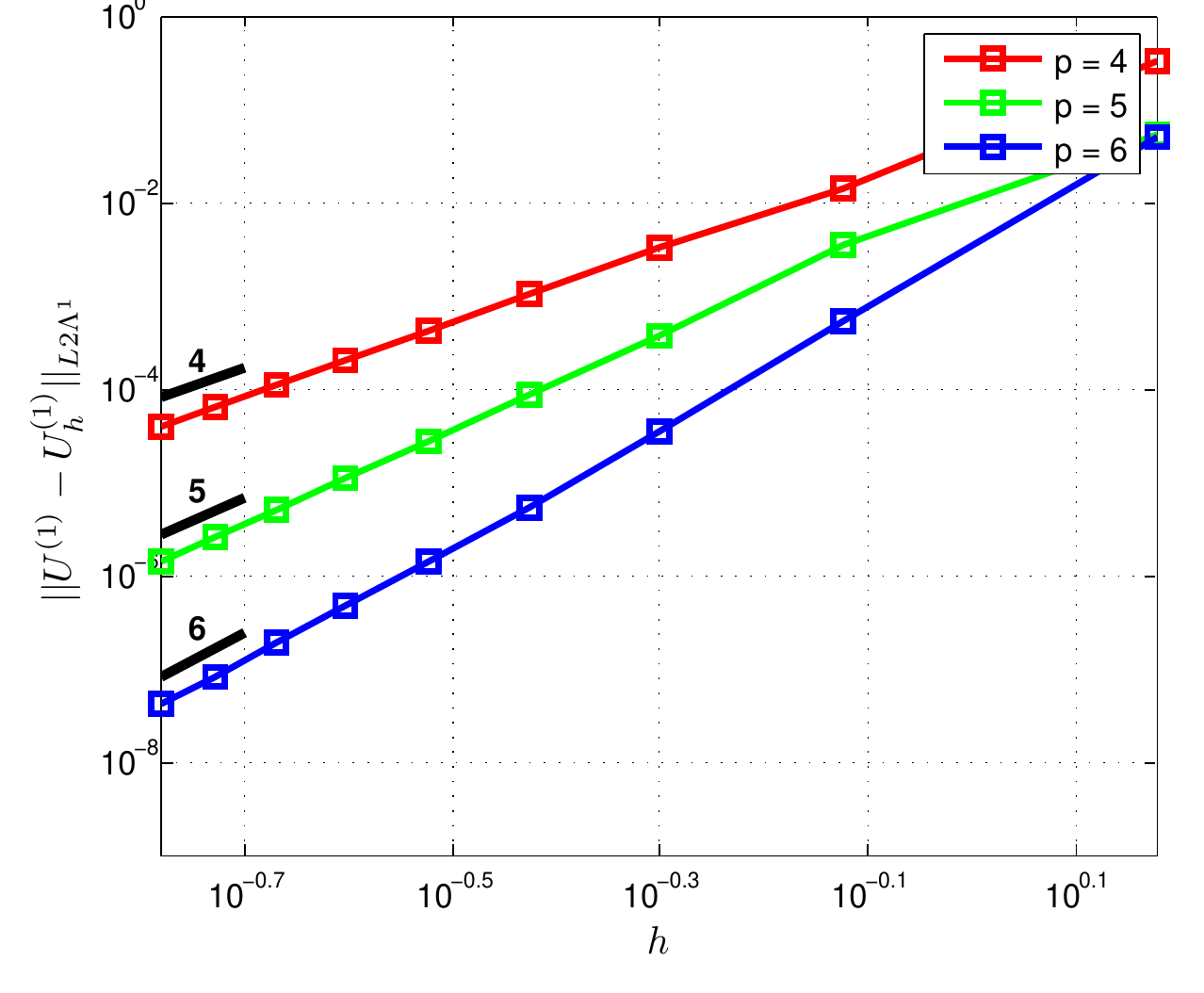}\label{fig::kovConvc}} 
  \subfigure[{Velocity, $p$-refinement}]{\includegraphics[trim = 10mm 4mm 14mm 9mm,clip,width=0.49\linewidth]{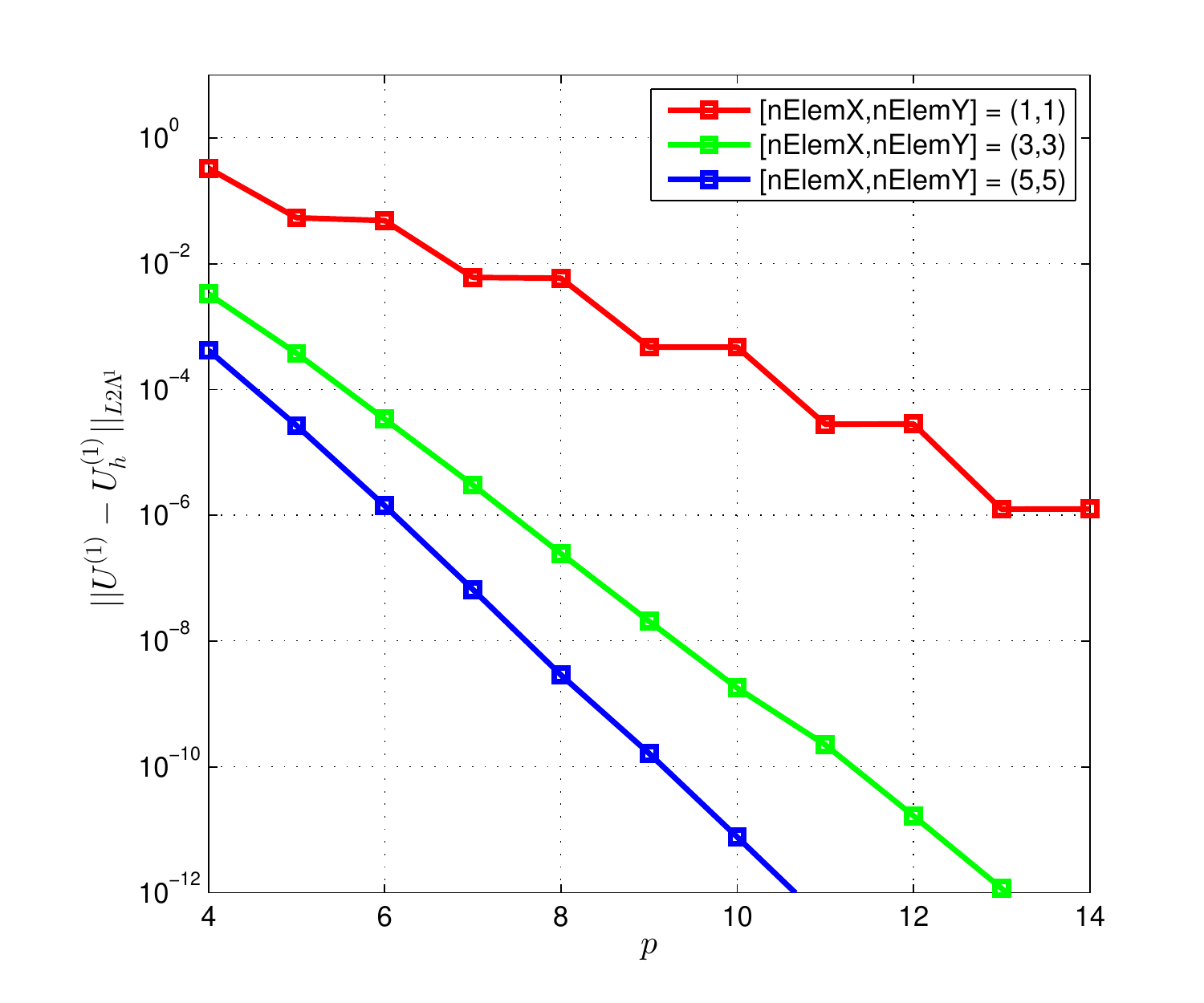}\label{fig::kovConvd}}
  \caption{Convergence plots for Kovasznay flow with mesh size, $h$, and order, $p$. Optimal rates are shown ({\color{black}{\solidrule}}) for the $h$-refinement cases. $nElemX$ and $nElemY$ refer to the number of elements in $X$ and $Y$ directions, and $p$ refers to the order of elements used.}
\label{fig::kovConv}
\end{figure}
In the above, $\beta^{(1)}$ is an arbitrary 1-form; $\tilde{M}_{11}$ and $\tilde{M}_{22}$ are mass-matrices for 1- and 2-forms on the staggered mesh; $C_{\vec{v}}$  is the convection matrix and depends on $\vec{v}$ (which can be retrieved from reconstruction of $\bar{v}$ using the edge basis, \citep{marcI}); $P^p_{\vec{w}}$ is the matrix which converts the scalar $\bar{p}$ to pressure-force 1-forms; $\tilde{B}_P$ and $\tilde{B}_T$ are the boundary integrals for pressure and stress, respectively, obtained from integration by parts; and ${D}_{21}$ and $\tilde{D}_{21}$ are incidence matrices which discretely represent the exterior-derivative with entries containing only $\{-1,0,1\}$.
The algebraic system thus obtained is solved for $\bar{v}$ and $\bar{p}$ for $\vec{w} = \{\vB{x}, \vB{y}\}$.
\section{Results}
\subsection{Kovasznay Flow}
Kovasznay flow is an analytical solution to Navier-Stokes' equations. The solution is $u = 1-e^{\lambda x}cos(2\pi y)$, $v = \frac{\lambda}{2\pi} e^{\lambda x}sin(2\pi y)$ and $p = \frac{1}{2} (1 - e^{2\lambda x})$, where $\lambda = \frac{1}{2\nu} - \sqrt{\frac{1}{4\nu^2} + 4\pi^2}$. The kinematic-viscosity chosen for this flow was $\nu = \frac{\mu}{\rho} = \frac{1}{40}$ and the computational domain considered was $\Omega = [-0.5 \; 1]\times[-0.5 \; -0.5]$.
\begin{figure}[htb!]
      \centering
    \subfigure[{Stream-function contours}]{\includegraphics[trim = 10mm 4mm 14mm 2mm,clip,width=0.49\linewidth]{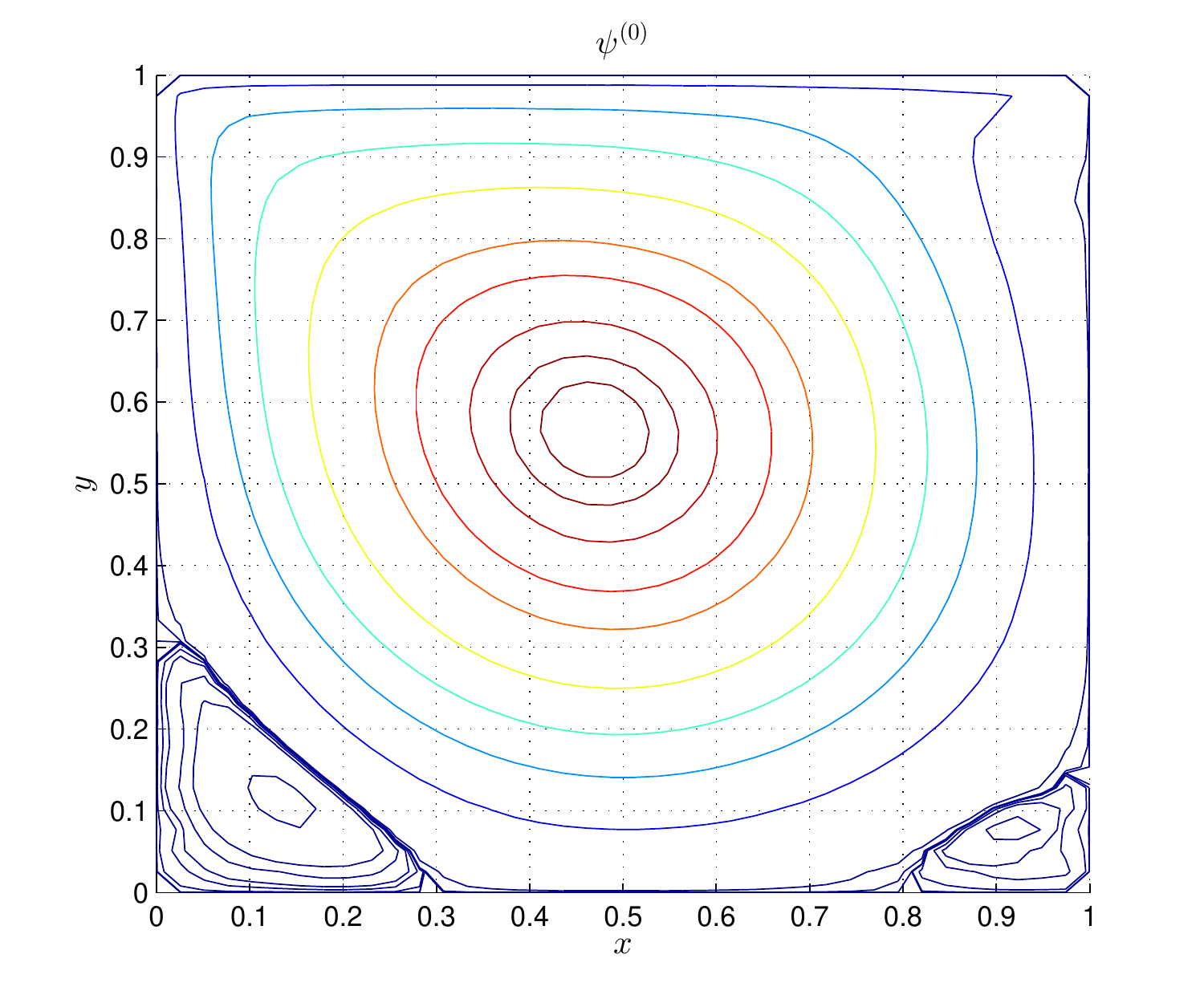}}
    \subfigure[{Pressure contours}]{\includegraphics[trim = 10mm 4mm 14mm 2mm,width=0.49\linewidth]{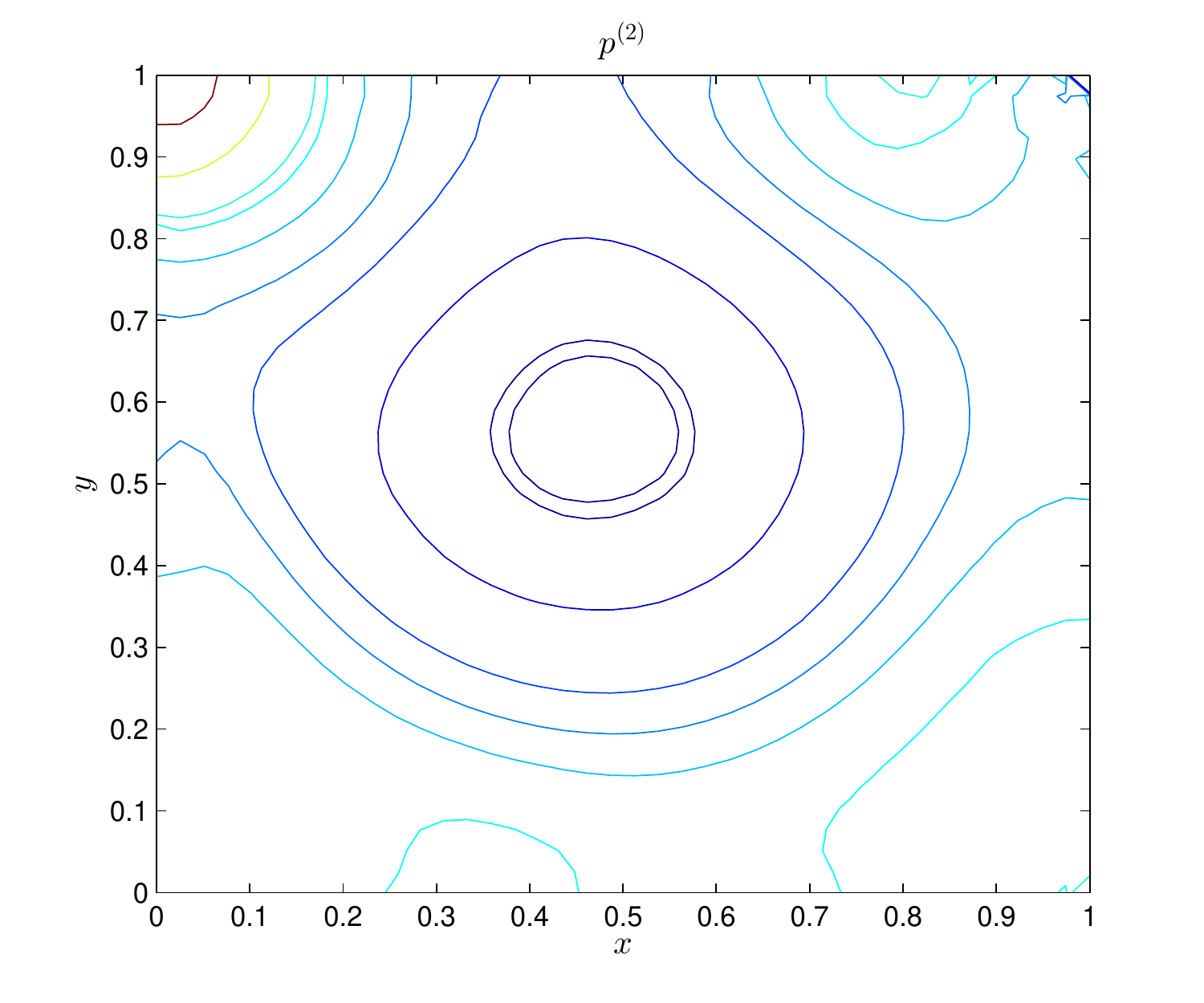}} \\
    \subfigure[{X-velocity at x = 0.5}]{\includegraphics[trim = 10mm 4mm 14mm 2mm,width=0.49\linewidth]{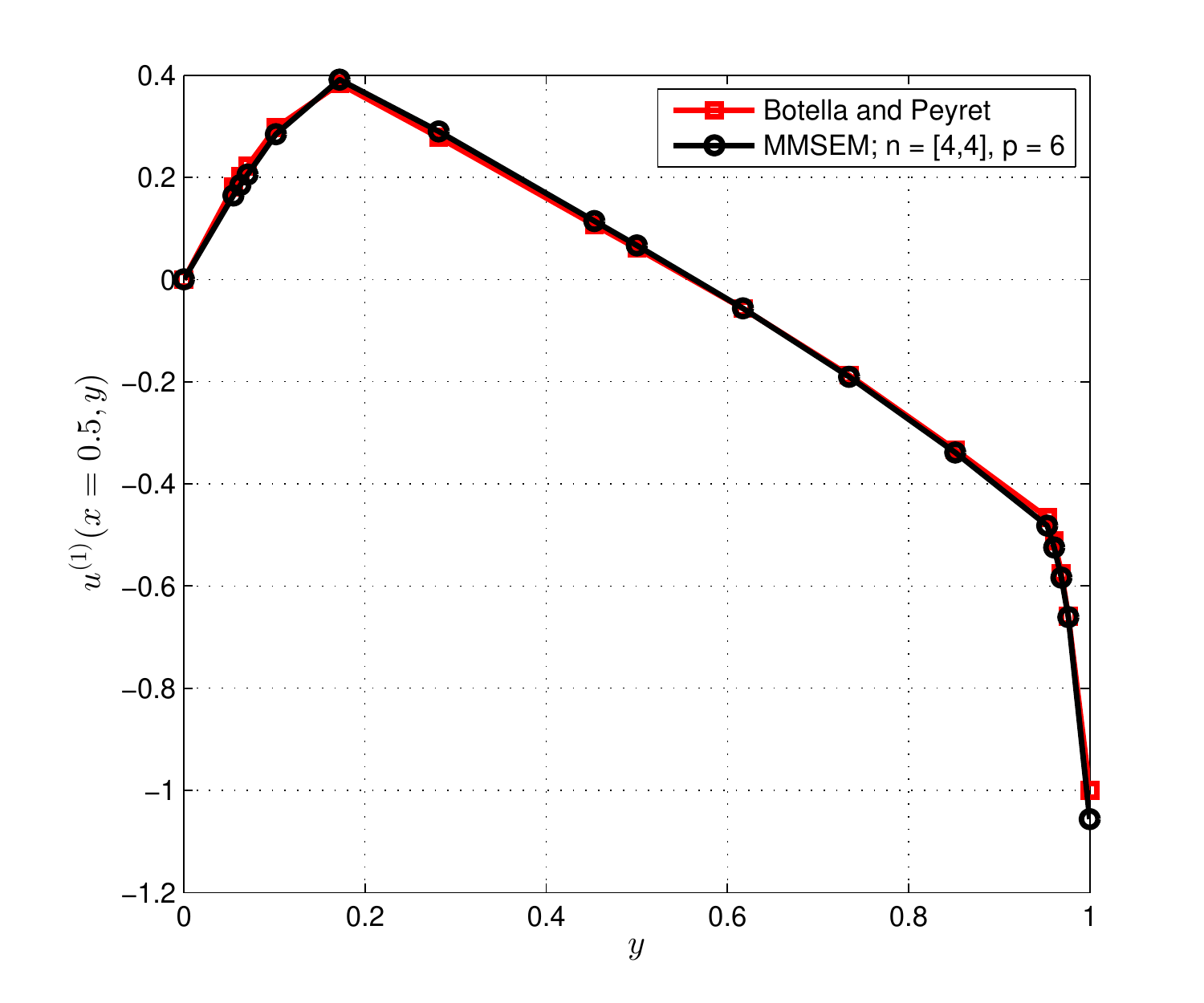}}
    \subfigure[{Y-velocity at y = 0.5}]{\includegraphics[trim = 10mm 4mm 14mm 2mm,width=0.49\linewidth]{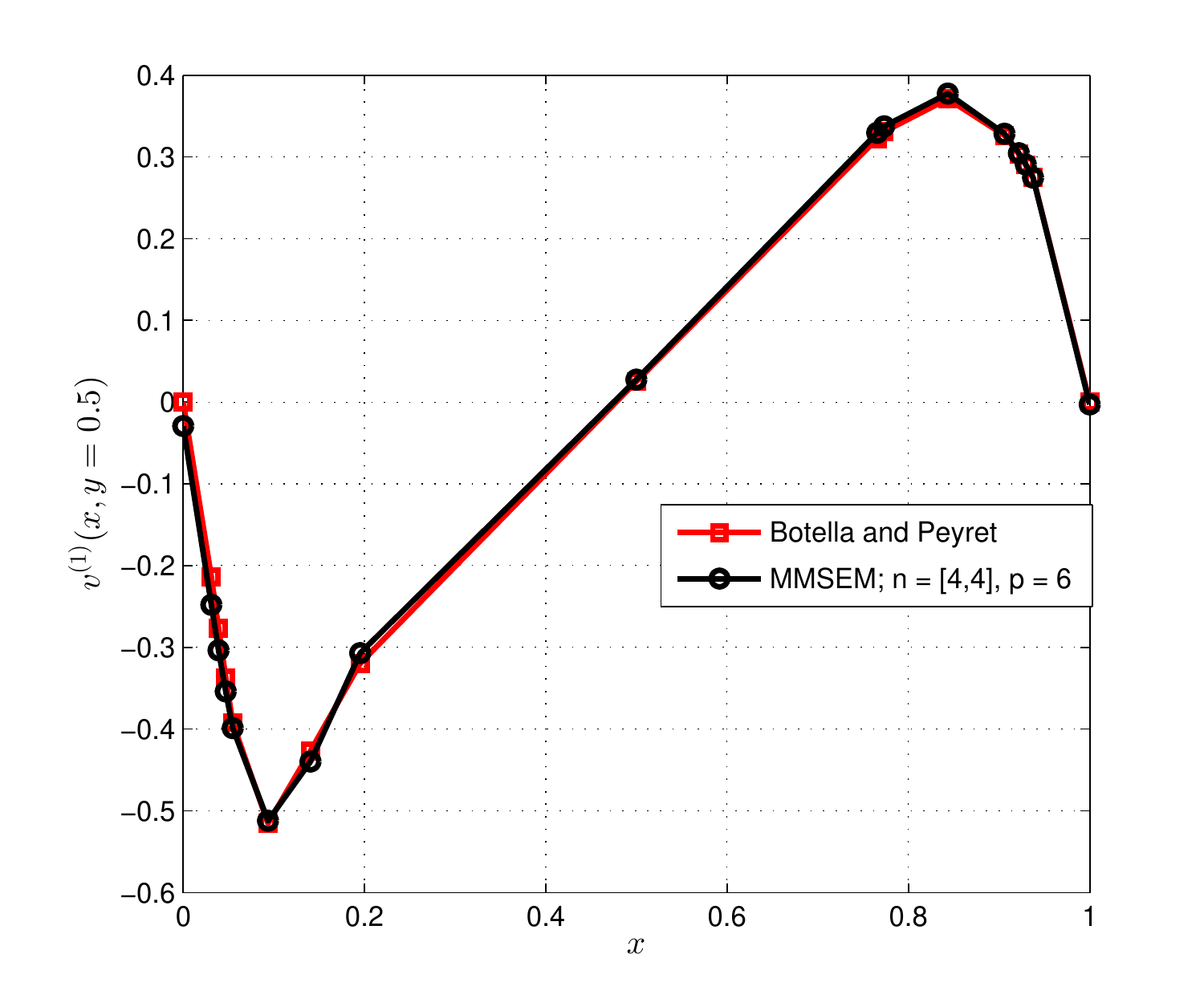}}
  \caption[Centerline velocities plotted for lid-driven cavity]{(Top) Streamfunction and Pressure contours with a single spectral element of order 16. (Bottom) Centerline velocities are plotted ({\color{black}{\solidrule}}) and compared with the solutions of \citep{Botella} ({\color{red}{\solidrule}}), and the solutions are found to be reasonably close. Mesh size is $4\times 4$ and made up of elements of order 6.}
\label{fig::lid}
\end{figure}    
The $h,p$-adaptivity plots for this problem are given in \figref{fig::kovConv} for pressures (\figref{fig::kovConva} and \figref{fig::kovConvb}) and velocities (\figref{fig::kovConvc} and \figref{fig::kovConvd}). It can be seen that the solutions converge exponentially and optimally. There is some stagnation observed in convergence for a mesh with a single element, and this is attributed to the fact that our basis may not be capturing certain modes (even/odd).

\subsection{Lid-driven Cavity Flow}
The second numerical test-case chosen was the classic Lid-driven cavity flow on a unit square domain with the top-lid velocity, $u_{L} = -1$ and a Reynolds number of 1000. The solutions for the pressure and streamfunction contours calculated for a single spectral element of order $p = 16$ are shown in the top-half of \figref{fig::lid}. Centerline-velocity solutions with a lower order of $p = 6$ but with multiple elements (4 x 4 mesh) and comparisons with the results of \citep{Botella} are also shown in the bottom-half of \figref{fig::lid}. Good agreement is seen between the benchmark results and our results.

\bibliographystyle{plainnat}
\bibliography{MyBib1}

\end{document}